\documentclass{amsproc}
\usepackage{euscript}
\usepackage{cases}
\usepackage{mathrsfs}
\usepackage{bbm}
\usepackage{amssymb}
\usepackage{amsfonts,amsmath,amsxtra,mathdots,mathabx}
\usepackage{color}
\usepackage{hyperref}
\usepackage{tikz}
\usepackage{appendix,upgreek}

\allowdisplaybreaks

\DeclareFontFamily{U}{matha}{\hyphenchar\font45}
\DeclareFontShape{U}{matha}{m}{n}{
	<5> <6> <7> <8> <9> <10> gen * matha
	<10.95> matha10 <12> <14.4> <17.28> <20.74> <24.88> matha12
}{}
\DeclareSymbolFont{matha}{U}{matha}{m}{n}

\DeclareMathSymbol{\Lt}{3}{matha}{"CE}
\DeclareMathSymbol{\Gt}{3}{matha}{"CF}

\def\valpha{\text{\scalebox{0.88}[1.02]{$\alpha$}}}   
\def\vepsilon{\upvarepsilon}

\def\vchi{\text{\raisebox{0.6 \depth}{\scalebox{0.9}[1.1]{$\chi$}}}} 
\def\vkappa{\text{{\scalebox{0.86}[1.1]{$\kappa$}}}} 
\def\vnu{\text{{\scalebox{0.9}[1]{$\nu$}}}} 
\def\uppii{\text{\scalebox{0.8}[0.96]{$\uppi$}}}

\newcommand{\BC}{{\mathbb {C}}}
\newcommand{\BQ}{{\mathbb {Q}}} 
\newcommand{\BR}{{\mathbb {R}}} 
\newcommand{\BZ}{{\mathbb {Z}}}

\newcommand{\GL}{{\mathrm {GL}}}
\newcommand{\PGL}{{\mathrm {PGL}}}
\newcommand{\SL}{{\mathrm {SL}}}

\newcommand{\ra}{\rightarrow} 
\def\viint{	\int \hskip -5 pt \int}

\def\sumx{\sideset{}{^\star}\sum}

\def\mod{\mathrm{mod}\,  }

\def\vwedge{\hskip -1.5pt \wedge \hskip -1.5pt}

\def\BSA{  \boldsymbol{\EuScript{A}} }
\def\HA{ \,\widehat{\phantom{A} } \hskip -9.5pt \BSA  }
\def\SHA{\, \widehat{\phantom{A} } \hskip -7.5pt \BSA}

\def\nd{\mathrm{d}}

\def\trh{ \mathrm{trh}}

\newcommand{\RF}{{\mathrm {F}}}

\def\SC{\text{\raisebox{- 2 \depth}{\scalebox{1.1}{$ \text{\usefont{U}{BOONDOX-calo}{m}{n}C} \hskip 0.5pt $}}}}

\def\SE{\text{\raisebox{- 2 \depth}{\scalebox{1.1}{$ \text{\usefont{U}{BOONDOX-calo}{m}{n}E} \hskip 0.5pt $}}}}

\def\SM{\text{\raisebox{- 2 \depth}{\scalebox{1.1}{$ \text{\usefont{U}{BOONDOX-calo}{m}{n}M} \hskip 0.5pt $}}}}

\def\lp {\left (}
\def\rp {\right )}

\def\boldJ {\boldsymbol J}

\renewcommand{\Re}{{\mathrm{Re} }}

\def\shskip{\hskip 0.5 pt}

\def\CaloO {\text{\raisebox{- 2 \depth}{\scalebox{1.1}{$ \text{\usefont{U}{BOONDOX-calo}{m}{n}O}  $}}}}

\def\SDH{\text{\raisebox{- 6 \depth}{\scalebox{1.06}{$ \text{\usefont{U}{dutchcal}{m}{n}H}  $}}}}

\newcommand{\delete}[1]{}

\theoremstyle{plain}

\newtheorem{thm}{Theorem} 
\newtheorem{lem}[thm]{Lemma}  \newtheorem{prop}[thm]{Proposition}
\newtheorem*{conj}{Conjecture}

\theoremstyle{remark} 
\newtheorem{remark}{Remark}[section]

\numberwithin{equation}{section}

\begin{document}

	\title[Spectral Large Sieve  for  $\mathrm{PGL}_2 (\mathbb{Z} {[i]})   \backslash \mathrm{PGL}_2 (\mathbb{C}) $]{{On the Spectral Large Sieve Inequality for $\mathrm{PGL}_2 (\mathbb{Z}[i]) \backslash \mathrm{PGL}_2 (\mathbb{C})$}}

\begin{abstract}
	  In this paper, we introduce concise expressions for the complex Bessel integral that enables us to improve the spectral large sieve inequality of Watt for $\mathrm{PGL}_2 (\mathbb{Z}[i]) \backslash \mathrm{PGL}_2 (\mathbb{C})$. Our result is optimistic in the case when the average of spectral parameters is on a square centered at $(0, 0)$. 
\end{abstract}

\author{Zhi Qi}
\address{School of Mathematical Sciences\\ Zhejiang University\\Hangzhou, 310027\\China}
\email{zhi.qi@zju.edu.cn}

\thanks{The author was supported by National Key R\&D Program of China No. 2022YFA1005300 and National Natural Science Foundation of China No. 12071420.}


\maketitle

\section{Introduction}

In two impressively long papers \cite{Watt-2-Large-Sieve, Watt-1-Large-Sieve}, Watt established several  spectral and Kloosterman  large sieve inequalities for Hecke congruence subgroups of $ \mathrm{PSL}_2 (\BZ [i]) $, and thus successfully extended the seminal work of Deshouillers and Iwaniec  \cite{D-I-Kuz} from $\BQ$ onto $\BQ (i)$. Among the various  problems in analytic number theory related to Kloosterman sums that  the results of \cite{D-I-Kuz} were applied to (see \cite[\S 1.5]{D-I-Kuz}), in his third long paper  \cite{Watt-3-Hecke},  Watt aimed at the mean-value problem of the (Dirichlet-polynomial) weighted {\it fourth} moment of the Riemann $\zeta$ function as in \cite{D-I-Mean-Riemann-1,Watt-4th-Riemann} (see also the earlier work of Iwaniec \cite{Iwaniec-Mean-Riemann}),  and obtained parallel results for Hecke $\zeta$ functions over  $\BQ (i)$. 

However, Watt indicated that his spectral large sieve inequalities in \cite[Theorem 1]{Watt-1-Large-Sieve} are `{not} quite a perfect analogy' of those in \cite[Theorem 2]{D-I-Kuz} and `so seem open to further improvement'. 

The purpose of this paper is to introduce  two concise expressions for the Bessel integral 
to simplify the analysis and to improve the results of Watt.  In particular, we obtain the true analog of  \cite[Theorem 2]{D-I-Kuz}  in the case when the average of   spectral parameters is on a  (large) square  centered at $(0, 0)$. As for the quadratic forms with Kloosterman sums,  our Proposition \ref{prop: quad form}    is   much  closer to     \cite[Proposition 3]{D-I-Kuz}   than \cite[Proposition 2]{Watt-1-Large-Sieve}. 

The motive of this paper is  a potential application  to the spectral {\it fourth} moment of $L$-functions, 
via a certain $\mathrm{PGL}_2 (\mathbb{Z}[i]) $ analog  of the functional equation (usually named as `spectral reciprocity' nowadays) of Kuznetsov (see \cite{Kuznetsov-4th-Moment,Ivic-Moment-1,Motohashi-Spectral-4th-Moment}),   yet to be established. 

As the main tool used in \cite{Watt-2-Large-Sieve, Watt-1-Large-Sieve}---the Kuznetsov formulae of Lokvenec-Guleska \cite{B-Mo2}---is established on an arbitrary imaginary quadratic field $\mathrm{F} = \BQ (\sqrt{-D})$,  one may extend the results of this paper to $\RF$ without any difficulty, at least  when the class number $h_{\RF} = 1$.

\subsection{Notation} \label{sec: notation}

As the notation in \cite{Watt-1-Large-Sieve} is quite heavy, while,  instead of Kloosterman sums,  our focus will be on analysis and subsequently on $L$-functions, so we shall follow closely the notation and setting in \cite{Qi-Hankel}, and, in particular, we shall work on $\mathrm{PGL}_2 (\mathbb{Z}[i]) \backslash \mathrm{PGL}_2 (\mathbb{C})$. The experienced readers should have no difficulty in the generalization to the case of Hecke congruence subgroups.

Denote $\RF = \BQ (i)$ and $ \CaloO  = \BZ [i]$. Let $\CaloO^{\times}$ be the group of units and $\CaloO   / \CaloO^{\times}$ denote the set of    ideals, usually identified with the Gaussian integers of argument in $[0, \pi/2)$. The letters $m$, $n$, and $c$ will be reserved for  (non-zero)  Gaussian integers in $\CaloO    $. Let $(n)$ denote the ideal $n \CaloO$.

Let $e [z] = \exp (2\pi i\Re (z))$.  For $  m,  n \in \CaloO $ and $c \in \CaloO \smallsetminus \{0\}$,      define the Kloosterman sum 
\begin{align}\label{1eq: defn Kloosterman}
	S   (m, n ; c ) = \sum_{ \valpha \delta \equiv 1 (\mathrm{mod} \, c) } e \bigg[  \frac {  \valpha m +   \delta n } {c} \bigg] .
\end{align} 
 


The   (unitary) characters on $\BC^{\times} $ are of the form     \begin{align}
	\vchi_{i \vkappa,  p} (z) = |z|^{  i \vkappa  } (z/|z|)^{  p}  , 
\end{align}  for   $ \vkappa   $ real and   $ p$ integral, 
so we introduce  $\BSA     = \BR \times  \BZ$ to parameterize the (Mellin) unitary dual of $\BC^{\times}  $.  
Define $ \HA = \BR \times \BR / \pi \BZ$, which will be considered as the (Fourier) dual of $  \BSA $. 
Let $\nd \mu (\vkappa,  p) $ or $\nd \widehat{\mu} (r, \omega)$  denote the usual Lebesgue measure on  $\BSA    $ or $ \HA$ respectively.  For simplicity, write $ \vchi_{  p} (z) = \vchi_{0,  p} (z)$.

For $p \in \BZ $,  define 
\begin{align}\label{1eq: zeta}
	\zeta (s, p) =  \sum_{(n) \shskip \subset \CaloO } \frac {\vchi_{4p} (n) } {|n|^{2s}}, \qquad \text{($\Re (s) > 1$)}, 
\end{align}
to be the Hecke $\zeta$ function associated to the Gr\"ossencharakter $\vchi_{4p} : (n) \ra (n/|n|)^{4 p}$, and   define 
\begin{align}\label{1eq: tau s}
	\tau_{s, p} (n) = \sum_{ (a) (b) = (n)} \vchi_{4p} (a/b) |a/b|^{2s} .
\end{align}  
Note that $\tau_{- s, - p} (n) = \tau_{s, p} (n)$. 

For $(\vkappa, p) \in \BSA$, let $\uppii_{i \vkappa, p} $ be the unitary principle series of $\mathrm{PGL}_2 (\BC)$: the unique infinite dimensional constituent of the representation   unitarily induced from the character
\begin{align*}
\vchi_{i \vkappa,  p}: 	\begin{pmatrix}
		\sqrt{z} & v \\
		& 1/\sqrt{z}
	\end{pmatrix}   \ra 
	|z|^{i \vkappa} (z/|z|)^{  p}. 
\end{align*}
Note that  $\uppii_{-i \vkappa, -p} \approx   \uppii_{i \vkappa, p} $.
As the Selberg conjecture holds for $\mathrm{PGL}_2 (\CaloO)$, we do not consider complementary series. 

Let  $  L^2_{c} (\mathrm{PGL}_2 (\CaloO) \backslash \mathrm{PGL}_2 (\BC))  $ be the space of  cusp forms for $ \mathrm{PGL}_2 (\CaloO) $.  Let $\Pi_{c} $ denote the discrete   spectrum of the irreducible constituents of $L^2_{c} (\mathrm{PGL}_2 (\CaloO) \backslash \mathrm{PGL}_2 (\BC))$. 

It may be assumed that each  $\uppii \in \Pi_{c}$ is Hecke invariant. Let $\lambda_{\uppii} (n)$ be the Hecke eigenvalues of $ \uppii $. It is known that $\lambda_{\uppii} (n)$ are all real.  Note that the $\mathrm{PGL}_2 (\CaloO) $-invariance ensures that $ \lambda_{\uppii} (n) $ only depends on the ideal $(n)$.

For $\uppii \in \Pi_{c}$, let $  (\vkappa_{\uppii}, p_{\uppii})$ be the parameter $(\vkappa, p) \in \BSA$ such that  $ \uppii \simeq \uppii_{i \vkappa , p }$.  

Let $L(s, \uppii)$ and $L (s, \mathrm{Sym}^2 \uppii)$ be the standard and the symmetric square $L$-functions attached to $\uppii$.   Note that  the root number of $ \uppii $ is equal to $(-1)^{ p_{\uppii}}$. 

Finally, define 
\begin{align*}
	\Pi_{c} (K, P) \hskip -1pt  = \hskip -1pt   \big\{ \uppii \in \Pi_{c} : |\vkappa_{\uppii} | \leqslant K,   |p_{\uppii} | \leqslant P \big\},  \quad  \BSA (K, P) \hskip -1pt     =     \hskip -1pt     \big\{ (\vkappa, p) \in \BSA : |\vkappa  | \leqslant K,   |p  | \leqslant P \big\}. 
\end{align*}


\subsection{Main Results}


\begin{thm}\label{thm: large sieve}
	Let $K, P , N  \geqslant  1/2$ {\rm(}say{\rm)} and $\vepsilon > 0$ be real numbers.  Let $\boldsymbol{a} : \CaloO / \CaloO^{\times}   \ra \BC$ be a sequence. Then both 
	\begin{align*}
\SC (K, P, N) =	& \sum_{\uppii \shskip \in \Pi_{c} (K, P) }  \frac { 1  } { L (1, \mathrm{Sym}^2 \uppii) } \bigg| \sum_{ N   <  |n|  \leqslant   \sqrt{2  } N }  a_{(n)} \lambda_{\uppii} (n)  \bigg|^2 , \\
\SE (K, P, N) =	&	\viint_{ \BSA (K/2, P/4) } \frac 1 { | \zeta (1+2 i \vkappa, 2 p)   |^2  } \bigg| \sum_{ N   <  |n|  \leqslant   \sqrt{2  }N }  a_{(n)} \tau_{i \vkappa, p} (n)  \bigg|^2  \nd \mu (\vkappa, p ) ,
	\end{align*}
	are majorized up to a constant depending on $\vepsilon$ at most by
	\begin{align}\label{1eq: bound, large sieve}
		\big(K^2 + P^2\big) \lp K P + K^{1+\vepsilon} +  P^{1+\vepsilon} + \frac{ N^{2+\vepsilon}} {KP}  \rp \|\boldsymbol{a}_{N} \|_2^2 \text{{\rm;}}
	\end{align}
here, and in what follows, we use the notation
\begin{align*}
	\| \boldsymbol{a}_N \|_2   = \bigg(  \sum_{  {N} <  |n|  \leqslant   \sqrt{2  }N } |a_{(n)}|^2 \bigg)^{1/2}. 
\end{align*}
\end{thm}

For comparison,  the corresponding bound of Watt \cite[Theorem 1]{Watt-1-Large-Sieve} in the case of $\mathrm{PSL}_2 (\CaloO)$ reads
\begin{align*}
	 \big(K^2 + P^2\big) \lp K P + \frac{ N^{2+\vepsilon}} {\sqrt{KP}}  \rp \|\boldsymbol{a}_{N} \|_2^2. 
\end{align*}

Note that the quality of the upper bound \eqref{1eq: bound, large sieve} decreases as $K/P + P/K$ increases. 
However, in the square case when $K = P$, we have the optimistic (conjectural) bound
\begin{align}\label{1eq: bound, K=P}
\sum_{\uppii \shskip \in \Pi_{c} (K, K) }  \frac { 1  } { L (1, \mathrm{Sym}^2 \uppii) } \bigg| \sum_{  {N} <  |n|  \leqslant   \sqrt{ 2 }N }  a_{(n)} \lambda_{\uppii} (n)  \bigg|^2 \Lt_{\vepsilon} 	\big(K^4 + N^{2+\vepsilon} \big) \|\boldsymbol{a}_{N} \|_2^2,
\end{align}
and this is the analog  of the large sieve inequality of Iwaniec for $\mathrm{PSL}_2 (\BZ)$ as in \cite[Theorem 1]{Iwaniec-Large-Sieve}. Moreover, if we had worked on a Hecke congruence subgroup of $\PGL_2 (\CaloO)$ or $\mathrm{PSL}_2 (\CaloO) $, then the $N^{2+\vepsilon}$ in \eqref{1eq: bound, K=P} is expected to be replaced by $|\mu (\mathfrak{a})|^2 N^{2+\vepsilon}$ if we average the $\boldsymbol{a}$-weighted Fourier coefficients at a cusp $\mathfrak{a}$  (for the details, see   \cite{Watt-1-Large-Sieve}); in this way, we  have indeed the analog  of \cite[Theorem 2]{D-I-Kuz}.

\subsection{Quadratic Forms with Kloosterman Sums}

\begin{prop}\label{prop: quad form}
	
	Let $\boldsymbol{b} : \CaloO / \CaloO^{\times}   \ra \BC$ be a sequence, $\theta \in \BC \smallsetminus \{0\} $, and $ c \in \CaloO \smallsetminus \{0\}$. Let $N > 0$. 
Write
	\begin{align*}
		B (\theta,   c, N) = {\mathop{\sum \sum}_{ {N} < |m| , |n|  \leqslant   \sqrt{  2 }N } } b_{(m)} \overline{b }_{(n)}   S (m, n; c)  \cos   \lp    4\pi  \Re  \lp \frac {\sqrt{m n}} {c} \theta   \rp    \rp .  
	\end{align*}
	For any $\vepsilon > 0$, we have 
	\begin{align}
		\label{1eq: Kloosterman, 1} &	   B (\theta,   c, N)   \Lt \tau^2 (c) |c| N^2 \|\boldsymbol{b}_N  \|_2^2, \\
		\label{1eq: Kloosterman, 2}&	   B (\theta,   c, N)   \Lt \big(|c|^2 + N^2 + |N c \shskip \theta|    \log (2 + |N \theta|    ) \big) \|\boldsymbol{b}_N \|_2^2, \\
		\label{1eq: Kloosterman, 3}&	  B (\theta,   c, N)    \Lt |N c/ \shskip \theta|   {N}^{ \vepsilon}  \|\boldsymbol{b}_N \|_2^2, 
	\end{align}
the last inequality holding provided that  $ |\theta | <   \sqrt[8]{2}$ and $  |c|  <  {N}$.  
\end{prop}

This may be considered as the true analog  of  \cite[Proposition 3]{D-I-Kuz}, 
 compared with  \cite[Proposition 2]{Watt-1-Large-Sieve}.

\subsection{Application: Bounding the Spectral  Moments of $L$-functions} 
Define the spectral $2q$-th moment of central $L$-values: 
		\begin{align*}
		\label{1eq: defn moment}	\SM_{2q} (K, P) =	  \sum_{\uppii \shskip \in \Pi_{c} (K, P) }  \frac { L (1/2,  \uppii)^{2q}  } { L (1, \mathrm{Sym}^2 \uppii) }  + \viint_{ \BSA (K/2, P/4) } \frac {|\zeta (1/2+i \vkappa, p) |^{4q}} { | \zeta (1+2 i \vkappa, 2 p)   |^2  }  \nd \mu (\vkappa, p )  .  
\end{align*}
\begin{thm}
	\label{thm: moments}
		Let $K, P   \geqslant  1/2$  and $\vepsilon > 0$ be real numbers. Then
		\begin{align}
			 &\SM_2 (K, P) \Lt K P \big(K^2+P^2\big)^{1+\vepsilon}, \qquad \SM_4 (K, P) \Lt   {\big(K^2+P^2\big)^{2+\vepsilon}}  ,  
		\end{align}
	with the implied constants depending on $\vepsilon$ only. 
\end{thm}

Thanks to the approximate functional equations, these bounds are a direct consequence of Theorem \ref{thm: large sieve} with $N $ up to $ (K^2+P^2)^{q/2+\vepsilon}$ ($q=1, 2$). 

\begin{conj}
	The mean Lindel\"of hypothesis for the {\it fourth} moment reads{\rm:} 
	\begin{align}
		 \SM_4 (K, P) \Lt   { K P \big(K^2+P^2\big)^{1+\vepsilon}}  .   
	\end{align}
\end{conj}

\subsection{Discussion on the Spherical Case} 
Let $\Pi_{c}^0$ denote the set of spherical $\uppii$ with $p_{\uppii} = 0$. By convention, let us identify $\Pi_{c}^0$ with an orthonormal basis of Hecke--Maass cusp forms $u_j $ on the Picard manifold $ \mathrm{PGL}_2 (\CaloO) \backslash \mathbb{H}^3  $. 
Let $\lambda_j = 1 + \vkappa_j^2$ ($\vkappa_j > 0$) and  $\lambda_{j} (n)$  denote the Laplace eigenvalue and the  the Hecke eigenvalues of $u_j$ respectively.  
Define
\begin{align*}
	\Pi_{c}^{0} (K ) = \big\{ u_j \in \Pi_{c}^0 : \vkappa_{j}   \leqslant K  \big\}. 
\end{align*}
On choosing $P = 1/2$ in Theorems \ref{thm: large sieve} and \ref{thm: moments}, we have
\begin{align}\label{1eq: bound, spherical}
	\sum_{  \vkappa_j \leqslant K }  \frac { 1  } { L (1, \mathrm{Sym}^2 u_j ) } \bigg| \sum_{  {N} <  |n|  \leqslant   \sqrt{ 2 }N }  a_{(n)} \lambda_{j} (n)  \bigg|^2 \Lt_{\vepsilon} 	\big(K^{3+\vepsilon} + K N^{2+\vepsilon} \big) \|\boldsymbol{a}_{N} \|_2^2,
\end{align}
\begin{align}\label{1eq: bounds L2, L4}
	\sum_{  \vkappa_j \leqslant K }  \frac { L (1/2,  u_j)^{2}  } { L (1, \mathrm{Sym}^2 u_j) } \Lt_{\vepsilon} K^{3+\vepsilon}, \qquad \sum_{  \vkappa_j \leqslant K }  \frac { L (1/2,  u_j )^{4}  } { L (1, \mathrm{Sym}^2 u_j) } \Lt_{\vepsilon} K^{4+\vepsilon}. 
\end{align}
However, from \eqref{1eq: bound, spherical}, with $N$ up to $K^{2+\vepsilon}$,  one may only get  $ O_{\vepsilon} (K^{5+\vepsilon})$ for the {\it fourth} moment of $L (1/2,  f)$. Actually, if one worked solely on  $ \mathrm{PGL}_2 (\CaloO) \backslash \mathbb{H}^3  $, in order to achieve the $K^{4+\vepsilon}$ in the second bound in \eqref{1eq: bounds L2, L4}, one would require the conjectural large sieve inequality:  
\begin{align*}
	\sum_{  \vkappa_j \leqslant K }  \frac { 1  } { L (1, \mathrm{Sym}^2 u_j ) } \bigg| \sum_{  {N} <  |n|  \leqslant   \sqrt{ 2 }N }  a_{(n)} \lambda_{j} (n)  \bigg|^2 \Lt_{\vepsilon} 	\big(K^3 +   N^{2+\vepsilon} \big) \|\boldsymbol{a}_{N} \|_2^2, 
\end{align*}
which seems to be out of reach at present, or, alternatively, the Vorono\"i summation for $a_{(n)} = \tau (n)$ and probably more (like `spectral reciprocity'). 
On the other hand, observe that the bound $K^{4+\vepsilon}$ is obtained  from the embedding $ \Pi_{c}^0 (K) \hookrightarrow \Pi_{c} (K, K) $ and the large sieve inequality in \eqref{1eq: bound, K=P}. 

\section{Proof of Proposition \ref{prop: quad form}} 

Our proof of Proposition \ref{prop: quad form} is essentially in parallel with that of  Proposition 3 in Deshouillers--Iwaniec \cite{D-I-Kuz}.  The readers who are familiar with their work may skip this section safely.  

For  $\theta \in \BC \smallsetminus \{0\}$, we investigate quadratic forms of the type
\begin{align}\label{2eq: defn of B}
	B (\theta, c, N) = {\mathop{\sum \sum}_{ {N} < |m| , |n|  \leqslant   \sqrt{2 }N} } b_{(m)} \overline{b }_{(n)}   S (m, n; c)   \cos \lp   4\pi  \Re \lp \frac {\sqrt{m n}} {c} \theta \rp \rp  .
\end{align}
Assume as we may that $ b_{(n)} = 0$   unless $  {N} <   |n|  \leqslant   \sqrt{ 2 }N  $. 

The first inequality \eqref{1eq: Kloosterman, 1}  follows easily from the Weil  bound:
  \begin{align*}
  	S (m, n; c)   \Lt \tau (c) {|(m, n, c)|} {|c|}. 
  \end{align*} 

The second inequality  \eqref{1eq: Kloosterman, 2} is a consequence of the   hybrid large sieve inequality: 
\begin{equation} \label{2eq: hybrid} 
		\viint_{ \BSA (T, M) } \sumx_{\delta (\mod c)}     \bigg|   \sum_{ n }  b_{(n)} \vchi_{i\vkappa,  p} (n) e \bigg[  \frac {\delta n} {c} \bigg] \bigg|^2 \nd \mu (\vkappa, p) \Lt  \big( (1+M) T |c|^2 +  {N^2}\big)   \| \boldsymbol{b}_N \|_2^2; 
\end{equation}
this  is the special case of \cite[Lemma 3.6]{Watt-1-Large-Sieve} with $\valpha = 1/2\pi$ and $\beta = -1$. To separate the variables $m$ and $n$ in $\sqrt{m n} \theta / c$, we need to use the Mellin technique. To this end, we invoke some results from \cite{Qi-GL(3)}. 
Let $R \Gt 1$. Then, for $R < |z| \leqslant \sqrt{2} R  $, it is proven in  \cite[Lemma 12.4]{Qi-GL(3)} that 
\begin{align}\label{2eq: Mellin}
	e [z] = \viint_{\BSA}  \xi_{R} (\vkappa, p) \overline{\vchi_{i\vkappa,  p} (z)} \nd \mu (\vkappa, p) , 
\end{align}
with 
\begin{equation*}
	\xi_{R} (\vkappa, p) \Lt \left\{ 
	\begin{split}
	&\displaystyle \frac {\log R} R, & & \text{ if } \sqrt{\vkappa^2 + p^2} \asymp R,  \\
	& (R + |\vkappa|+|p|)^{-A}, & & \text{ if otherwise.  }
	\end{split}
	\right.
\end{equation*}
For the (easier) case $R \Lt 1$, it follows from \cite[Lemma 12.3]{Qi-GL(3)} that \eqref{2eq: Mellin} is still valid with 
\begin{equation*}
	\xi_{R} (\vkappa, p) \Lt   ( 1 + |\vkappa|+|p|)^{-A}. 
\end{equation*}
In our setting, let $R =    |N \theta /c|$ and choose  $A = 3$ (say).  Hence, by \eqref{2eq: defn of B} and \eqref{2eq: Mellin} we infer that
\begin{align*} 
	B (\theta, c, N)   =    \viint_{\BSA}   \xi_{   |N \theta/ c|} (\vkappa, p) {\mathop{\sum \sum}_{m, n} }  b_{(m)} \overline{b }_{(n)}  \overline{\vchi}_{i\vkappa,  2 p} \lp \frac {\sqrt{m n}} {c} \theta \rp   S (m, n; c)  \nd\mu (\vkappa, p)   .
\end{align*}
By \eqref{1eq: defn Kloosterman}, \eqref{2eq: hybrid}, and Cauchy--Schwarz, along with the bounds for $\xi_{R} (\vkappa, p)$,  we have 
\begin{align*}
	B (\theta,   c, N)   \Lt \big(  |N c \shskip \theta|    \log   |N\theta|   + N^2 \big) \|\boldsymbol{b}_N \|_2^2,
\end{align*}
if $ |N \theta |    \Gt |c|$, and 
\begin{align*}
	B (\theta,   c, N)   \Lt \big(  |c|^2 + N^2 \big) \|\boldsymbol{b}_N \|_2^2,
\end{align*}
if  $ |N \theta |    \Lt |c|$, as desired. 

Now we proceed to prove \eqref{1eq: Kloosterman, 3} for $|c| < N^{1 -\vepsilon}$; the remaining range is already covered by \eqref{1eq: Kloosterman, 2}. Let $  \eta \in C_c^{\infty} [1/\sqrt[6]{2}, 2]$ with $\eta \equiv 1$ on $[1, \sqrt{2}]$. By  Cauchy--Schwarz, 
\begin{equation*}
	\begin{split}
		|B (\theta,   c, N)|^2 & \leqslant \|\boldsymbol{b}_N     \|_2^2 \sum_{\pm} \frac 1 2   \sum_{m} \eta    \bigg(   \frac {|m|} { {N}}   \bigg)  \bigg| \sum_{n} b_{(n)} S (m, n; c)  e \bigg[ \hskip -1pt \pm 2    \frac {\sqrt{m n}} {c} \theta  \bigg]  \bigg|^2 \\
	& = \hskip -1pt   16  \|\boldsymbol{b}_N \hskip -1pt  \|_2^2 \hskip -2pt  \mathop{\sumx\sumx}_{\delta_1    , \delta_2 (\mod c)} \hskip -1pt \mathop{\sum \sum}_{(     n_1  ), (   n_2  )} \hskip -2pt b_{(    n_1  )} \overline{b}_{(    n_2  )} e\bigg[ \hskip -1pt \frac {n_1 \delta_1 -   n_2 \delta_2} {c} \hskip -1pt \bigg] C (\delta_1,   \delta_2,   (   n_1   ),   (   n_2     )   )   , 
	\end{split}
\end{equation*}
where  $C = C (\delta_1, \delta_2, (n_1), (n_2))$ is the sum
\begin{align*}
	C = \sum_{m \shskip \in  \shskip \CaloO} f (m) =    \sum_{m \shskip \in  \shskip \CaloO} \eta    \bigg(    \frac {|m|} { {N}}   \bigg) e \Big[   \frac {\valpha_1 -\valpha_2} {c} m   \Big] \cos   \bigg(    4\pi \Re \bigg(   \frac { \sqrt{n_1}     -      \sqrt{n_2} } {c} \theta\sqrt{m} \bigg)   \bigg),
\end{align*}
with $\valpha_1 \delta_1 \equiv \valpha_2 \delta_2 \equiv 1 (\mod c)$ and $\arg (\sqrt{n_1})$, $\arg (\sqrt{n_2}) \in [0, \pi/4)$. 
For brevity, let us write 
\begin{align*}
	a = \frac {\valpha_1 -\valpha_2} {c} , \qquad b = \frac { \sqrt{n_1}    -     \sqrt{n_2} } {c} \theta . 
\end{align*}
For the sum $C$ we apply the Poisson summation over $\CaloO$, obtaining
\begin{align*}
C = \sum_{u \shskip \in  \shskip \CaloO } \hat{f} (u) =	& \sum_{u \shskip \in  \shskip \CaloO }   \viint_{\BC} \eta     \bigg(    \frac {|z|} { {N}}   \bigg) e  [   (a - u) z    ] \cos  (  4\pi \Re  ( b \sqrt{z}  ) ) \frac {i \nd z \vwedge \nd \widebar{z}} {2 }  , 
\end{align*}
and, on the change $\pm \sqrt{z} \ra \sqrt{N} z$, 
\begin{align*}
\hat{f} (u) =  {N^2}  \viint_{\BC} \eta   (|z|^2) |z|^2 \cdot  e   \big[   (a - u)  {N} z^2   -  2   b \sqrt{N} z  \big]   {i \nd z \vwedge \nd \widebar{z}}   .
\end{align*}
Note that 
\begin{align*}
\frac { {N} |z|} {|c|} -	\big|  b \sqrt{N} \big|  > \lp \frac 1 {\sqrt[12]{2}} -    \textstyle \sqrt{\sqrt{2} + 1 -   \sqrt[4]{8}} \cdot \sqrt[8]{2} \rp  \frac {  {N} } {|c| } \Gt  \frac {  {N}} {|c|} . 
\end{align*}
Thus for $u \neq a$, we have   
\begin{align*}
	 \big| (a-u)  {N} z - b \sqrt{N}\big| \Gt  |(a-u)  {N}  |.
\end{align*}
By \cite[Lemma 7.4]{Qi-GL(3)}, we infer that
\begin{align*}
	C = \hat{f} (a) + O \big( N^2 (|c|/ {N})^{A}  \big) = \hat{f} (a) + O \big( 1/|c|^2 \big),
\end{align*}
if we choose $A = [2/\vepsilon]$. It remains to consider 
$$\hat{f} (a) =  {N^2}  \viint_{\BC} \eta   (|z|^2) |z|^2 \cdot  e   \big[   2   b \sqrt{N} z  \big]   {\nd z \vwedge \nd \widebar{z}}  . $$  Since $a \in \CaloO$ in this case, it implies $ \valpha_1 \equiv \valpha_2 (\mod c)  $ and consequently  $ \delta_1 \equiv \delta_2 (\mod c)  $. As $\hat{f} (a)$ is independent on $\delta_1$, so we can first carry out the summation over $\delta_1$, and
\begin{align*}
	\bigg| \  \sumx_{\delta_1 (\mod c)} e \Big[ \frac {n_1 - n_2} {c} \delta_1 \Big] \bigg| = |S (n_1-n_2, 0; c)| \leqslant |(n_1-n_2, c)|^2. 
\end{align*} 
Besides the trivial bound $\hat{f} (a) \Lt N^2$, we also have
\begin{align*}
	\hat{f} (a) \Lt   { {N}}  / { |b|^2} \Lt \frac {|Nc|^2 } {|\theta|^2 |n_1-n_2|^2} ,
\end{align*} 
by partial integration, provided $n_1 \neq n_2$. Gathering all the above results together, we conclude that
\begin{align*}
		|B (\theta,   c, N)|^2 & \Lt |Nc|^2   \|\boldsymbol{b}_N \|_2^4 +   {|N c/\theta|^2  } \mathop{\sum\sum}_{(  n_1 )\neq ( n_2  ) }   \big|b_{(  n_1 )}  {b}_{( n_2 )} \big| \frac {|(n_1-n_2, c)|^2} {|n_1-n_2|^2}  \|\boldsymbol{b}_N \|_2^2 \\
		& \Lt {|N c/\theta|^2 } N^{ \vepsilon} \| \boldsymbol{b}_N \|_2^4. 
\end{align*}
Thus the proof of \eqref{1eq: Kloosterman, 3} and of Proposition \ref{prop: quad form} is completed. 

\section{\texorpdfstring{Spectral Kuznetsov Trace Formula for $\mathrm{PGL}_2 (\BZ [i])$}{Spectral Kuznetsov Trace Formula for PGL2(Z[i])}} 

Our main tool  is the spectral Kuznetsov trace formula of Bruggeman and Motohashi over the Gaussian field $\mathrm{F}$ \cite[Theorem 10.1]{B-Mo}; this was generalized to Hekce congruence groups for imaginary quadratic $\mathrm{F}$  \cite[Theorem 11.3.3]{B-Mo2}. 


For $\vnu \in \BC$, $  p \in \BZ  $, and $z \in \BC \smallsetminus \{0\}$,  define 
\begin{equation}\label{0def: J mu m (z)}
	J_{\vnu ,   p} (z) = J_{\vnu + p }    (z)   J_{\vnu -  p  }    ({\widebar z} ) ,
\end{equation} 
\begin{equation}\label{0eq: defn of Bessel}
	\boldJ_{ \vnu,   p} (z)  =   \frac {2\pi^2} {\sin (\pi \vnu)} \lp J_{-\vnu,\shskip  -p} (   z) - J_{\vnu,   p} (    z)  \rp.    
\end{equation}
The Bessel functions $J_{\vnu ,\shskip  p} (z)$ and  $\boldJ_{ \vnu,   p} (z)$ are even and real analytic in $z$, and (complex) analytic in $\vnu$.  
For the definitions above, the reader is referred to  \cite[\S   15.3]{Qi-Bessel}, \cite[(6.21), (7.21)]{B-Mo} or \cite[(4.58), (9.26)]{B-Mo2}.\footnotemark

\footnotetext{The Bessel functions for $\GL_2 (\BC)$ or $\SL_2 (\BC)$ are derived by different means in \cite{Qi-Bessel} and \cite{B-Mo,B-Mo2}.  
	It should be remarked that it is unnecessary and inconvenient to  modify   $$J_{\vnu} (z) = \sum_{n=0}^{\infty} \frac {(-)^n (z/2)^{\vnu +2n}} {n! \Gamma (\vnu + n + 1)} $$ by the factor $(z/2)^{-\vnu}$ (to make it analytic on $\BC$) as in \cite[(6.21)]{B-Mo}. This would destroy the Bessel equation!  }

For simplicity, define $\mathscr{H}$ to be the space of even functions $h (\vkappa, p)$  on $\BSA$   that admit  an entire analytic continuation $ h (\vkappa + i \sigma, p) $ so that it  decays rapidly in both $ \vkappa$ and $ p $,  uniformly for $\sigma$ on   bounded intervals.

\begin{lem}\label{lem: Kuznetsov}
Let the notation be as in {\rm \S \ref{sec: notation}}. 	Let $ h  \in \mathscr{H} $. 
	For $m , n \in \CaloO \smallsetminus \{0\}$, we have the identity{\rm:} 
	\begin{equation}\label{1eq: Kuznetsov} 
		\begin{aligned}
		  \sum_{\uppii \shskip \in \Pi_{c} }  &    \frac {  \lambda_{\uppii} ( m )     \lambda_{\uppii}  ( n )} {L(1, \mathrm{Sym^2} \uppii )}   h  ( \vkappa_{\uppii}, p_{\uppii} ) +   \frac 1 {\pi} \viint_{\BSA}  
			\frac {\tau_{i \vkappa, p} (m )  \tau_{i \vkappa, p} ( n )} {|\zeta (1+2i \vkappa, 2p)|^2} h ( 2\vkappa, 4 p )  \shskip   \nd  \mu (\vkappa, p)  \\
			& =        \frac {1} { 16\pi^3 }  \delta_{(m),   (n)} \cdot  \SDH  + \frac 1 {64 \pi^3 }  \sum_{ \epsilon \shskip  \in   \CaloO^{\times} \hskip -1pt / \CaloO^{\times 2} } \sum_{c  \shskip \in   \CaloO \smallsetminus \{0\} } \frac {S  (  m , \epsilon  n  ; c  ) } { |c|^2  } \SDH  \bigg( \frac { 2\pi \sqrt{\epsilon m n} } {    c     }   \bigg),
		\end{aligned}
	\end{equation}
	where  $\delta_{(m),   (n)}$ is the Kronecker $\delta$ symbol for ideals, $ \SDH  $ and $ \SDH (z)$ are the Plancherel and Bessel integrals defined by  
	\begin{align}\label{1eq: defn Bessel integral}
		\SDH  = \hskip -2pt  \viint_{\BSA}   h (\vkappa, p)  (\vkappa^2 + p^2 )  \nd  \mu (\vkappa, p), \quad \SDH (z) = \hskip -2pt  \viint_{\BSA}   h (\vkappa, p)  \boldsymbol{J}_{  i \vkappa, p} ( z )  (\vkappa^2 + p^2 )  \nd  \mu (\vkappa, p)  . 
	\end{align} 
\end{lem}

The reader might have noted some differences between \eqref{1eq: Kuznetsov} and \cite[(10.1)]{B-Mo}.  The  $\epsilon$-sum on the right of \eqref{1eq: Kuznetsov}  is due to the translation from $\mathrm{PSL}_2 (\CaloO)$ to $\PGL_2 (\CaloO)$ (see   \cite[(3.17)]{CI-Cubic} and \cite[Proposition 1]{Venkatesh-BeyondEndoscopy}). The harmonic weight  is expressed here as   (a multiple of) $ 1/  L(1, \mathrm{Sym^2} \uppii)$ (see \cite[Lemma 5]{SH-Liu-Maass})  in place of their $|c_{V} (1)|^2$ (their $V$ is the underlying space of our $\uppii$) by standard Rankin--Selberg calculations.

\section{Analysis for the Bessel Integral}


Let    
\begin{equation}\label{4eq: choice of h} 
	h( \vkappa,  p ) = \exp \lp - \lp \frac {\vkappa } {K} \rp^2 - \lp \frac {p } {P } \rp^2 \rp .
\end{equation}
This section is devoted to the study of the Bessel integral
\begin{equation}\label{5eq: Bessel H(z)}
	\SDH (z) = \hskip -2pt  \viint_{\BSA}   h (\vkappa, p)  \boldsymbol{J}_{  i \vkappa, p} ( z )  (\vkappa^2 + p^2 )  \nd  \mu (\vkappa, p)  .
\end{equation}

For $(\vkappa, p) \in  {\BSA}$, in the polar coordinates, we have the following integral representation of $\boldJ_{ i \vkappa,   p} (z)  $:
\begin{equation}\label{5eq: integral representation}
	\boldJ_{i \vkappa ,   p}  (  x     e^{ i \phi} )   = 4 \pi i^{2p}  \int_{0}^\infty  \overline{E  (ye^{i\phi})}{}^{ 2p}  J_{2p} \hskip -1pt \left(    {x} Y (y e^{i\phi}) \right) y^{2i\vkappa - 1}  \nd y,
\end{equation} 
with 
\begin{align}\label{5eq: Y and E}
	Y (z) = | {z}+1/  {z}|, \qquad E (z) = \frac { {z}+1/ {z}} {| {z}+1/ {z}|}. 
\end{align}
See \cite[Theorem 12.1]{B-Mo} and \cite[Corollary 6.17, Example 15.6]{Qi-Bessel}. 
The integral on the right of \eqref{5eq: integral representation} is absolutely convergent due to the $        1 /  {\sqrt{x}}  $ decay of $J_{2p} (x)$. 

Moreover, we have the   Bessel integral representation  (see \cite[\S 2.2]{Watson}):
\begin{align}\label{5eq: Bessel}
	J_{2p} ( x) = \frac 1 {2\pi i^{2p}} \int_0^{2\pi} \exp ( { 2i p \omega + i x \cos \omega} ) \nd \omega. 
\end{align} 
Thus $ |J_{2p} ( x)|  \leqslant 1$. Also,  it follows from \cite[\S 7.13.1]{Olver} that $ J_{2p} ( x) \Lt 1/\sqrt{x}$ provided that $ x > 1+p^2 $. Consequently,   we have uniformly 
\begin{align}\label{5eq: Bessel bound}
	J_{2p} (x) \Lt \frac {\sqrt{1+p^2}} {\sqrt{x}} . 
\end{align}

First we   use $y = \exp (r)$ to transform \eqref{5eq: integral representation} into
\begin{equation}\label{10eq: Bessel}
	\boldJ_{i \vkappa ,   p}  (  x     e^{ i \phi} )   = 4 \pi i^{2p}  \int_{-\infty}^\infty   \widebar{\vchi}_{2p}  (\cosh (r + i\phi) )   J_{2p} \hskip -1pt \left( 2  {x} |\cosh (r + i\phi)|  \right) \exp ({2 i r \vkappa})  \nd r . 
\end{equation}  
For  $|r| > 1$,  \eqref{5eq: Bessel bound} yields the uniform bound 
\begin{align*}
	J_{2p} \hskip -1pt \left( 2   {x} |\cosh (r + i\phi)| \right) \Lt \frac {\sqrt{1+p^2}} {\sqrt{x \cosh r}  } . 
\end{align*}
Thus the triple integral obtained by inserting \eqref{10eq: Bessel} into \eqref{5eq: Bessel H(z)} is absolutely convergent, and  hence it is permissible  to change the order of integration. 

For  $a \in \BC \smallsetminus \{0\}$, the Bessel formula \eqref{5eq: Bessel} yields
\begin{align}\label{10eq: Bessel, 2}
	\widebar{\vchi}_{2p} (a) J_{2p} (  |a|) = \frac 1 {2 \pi i^{2p}	} \int_0^{2\pi} \exp ( {2 i p \omega + i   \mathrm{Re} ( a e^{i \omega}  )} ) \nd \omega ,
\end{align}
on the change $ \omega \ra \omega + \arg (a)$. 

Now let the $r$-integral as in \eqref{10eq: Bessel} be the outermost, and apply \eqref{10eq: Bessel, 2} with  $ a = 2 x \cosh (r + i\phi)$ so that $ \mathrm{Re} ( a e^{i \omega}  ) = 2 x \shskip \trh (r, \omega; \phi)$, defined by
\begin{align}\label{5eq: trh function}
	\trh  (r, \omega; \phi) =    \cosh r \cos \omega \cos \phi - \sinh r \sin \omega \sin \phi .
\end{align}      It follows that the Bessel integral $\SDH	(  x     e^{ i \phi} )$ is equal to
\begin{align*}
	 4  \viint_{ {\SHA}  }  \cos ( 2 x \trh (r, \omega; \phi)) \lp \viint_{\BSA    }     \exp ( {2 i \vkappa r + 2 i p \omega} ) h (\vkappa, p)    (\vkappa^2 + p^2 )  \nd  \mu (\vkappa, p) \rp  \nd \widehat{\mu} (r,  \omega ) , 
\end{align*}
in the notation of \S \ref{sec: notation} ($\BSA = \BR \times \BZ$ and $\HA = \BR \times \BR /\pi \BZ$). Finally, by an evaluation of the inner integral, using 
\begin{align*}
	\int  \exp \big( 2 ir \vkappa       - (\vkappa /K)^2 \big)  \nd \vkappa = \sqrt{\pi} K \exp \big( \hskip -1pt    -  (K { r }  ) ^2   \big) ,
\end{align*}
and (by Poisson)
\begin{align*}
	 \sum   \exp \big(     2 i  \omega p   -  (  p / P  )^2 \big) = \sqrt{\pi} P \sum  \exp \big( \hskip -1pt   - (P (\omega + \pi q) )^2 \big) ,
\end{align*}
we conclude with the following integral representation for $\SDH (z)$.

\begin{lem}\label{prop: integral repn}
We have
	\begin{align}\label{4eq: integral H(z)}
		 \SDH (z) =  -  \viint_{ {\SHA}  }  \cos ( 2 \Re (z \trh (r, \omega)) )    \big(k''(r) \theta (\omega) + k (r) \theta'' (\omega)\big) \nd \widehat{\mu} (r,  \omega ) , 
	\end{align}
with
\begin{align}\label{10eq: trh function, 0}
	\trh  (r, \omega) =    \cosh r \cos \omega + i \sinh r \sin \omega,
\end{align}  
\begin{align}\label{4eq: defn of h}
	k  (r ) =  \sqrt{\pi} K \exp \big( \hskip -1pt    -  (K r )^2   \big), \qquad 	
	\theta (\omega) =   \sqrt{\pi} P \sum  \exp \big( \hskip -1pt   - (P (\omega + \pi q) )^2 \big)    . 
\end{align} 
\end{lem}

\begin{remark}
	The expression of $  \boldJ_{ i \vkappa ,   p} (z) $ given by {\rm(\ref{5eq: integral representation}, \ref{5eq: Y and E})} was used in a quite different way by Watt{\rm;} see Lemmas {\rm 5.4--5.6} in {\rm\cite{Watt-1-Large-Sieve}}. However, the use of this here is only for the sake of strictness: one may  obtain {\rm\eqref{4eq: integral H(z)}} more directly  by inserting into  \eqref{5eq: Bessel H(z)}  the {\it formal} integral representation\footnote{Actually, in the   setting of $\GL_n$, the author discovered \eqref{0eq: formal integral} prior to  (\ref{5eq: integral representation}, \ref{5eq: Y and E}) as in  \cite[\S \S  5, 6]{Qi-Bessel}.   This kind of {\it formal} integrals has played a central role   in the theory of {\it fundamental}  Bessel functions for $\GL_n (\BR) $ and $\GL_n (\BC)$; see  \cite[\S \S 7--9, 12]{Qi-Bessel}. }
	\begin{equation}\label{0eq: formal integral}
		\boldJ_{ i\vkappa, p }  (z)  =   \viint_{\BC^{\times}} \exp (i  \Re (z u + z / u)) \shskip 	\vchi_{i\vkappa,  p} (u^2 ) \frac {i \nd      u  \vwedge \nd      \widebar{u }} {|u|^2}; 
	\end{equation}
 this integral is {\it not} absolutely convergent for any value of $(\vkappa, p)$. 
	
\end{remark}

By partial integration, along with the identity
\begin{align*}
	\big(\hskip -0.5pt   (\partial / \partial r)^2 \hskip -1pt   + \hskip -1pt    (\partial / \partial \omega)^2 \hskip -0.5pt    \big) \hskip -1pt     \cos  ( 2 \Re (\hskip -0.5pt   z \trh (r, \omega) \hskip -0.5pt ) \hskip -0.5pt  ) \hskip -1pt     = \hskip -1pt      - 4    (\sinh^2 \hskip -1pt   r  \hskip -1pt      + \hskip -1pt      \sin^2 \hskip -1pt    \omega  )	|z|^2 \hskip -1pt   \cos  ( 2 \Re (\hskip -0.5pt   z \trh (r, \omega) \hskip -0.5pt   ) \hskip -0.5pt ), 
\end{align*} 
we deduce easily the following variant of Lemma \ref{prop: integral repn}. 

\delete{Moreover, by  repeating partial summation to $\theta (\omega)$ as defined in \eqref{4eq: defn of h}, we infer that $ h (r) \theta (\omega)$  is Schwartz of deviation $(1/K, 1/P)$ and magnitude $KP (K^2+P^2)  $   in the sense that 
\begin{align*}
	\frac 1 {K^i P^{j}} \frac {\partial^{i+j} h (r) \theta (\omega)  } {\partial r^i \partial \omega^j} \Lt_{i, j, A} KP   \lp \frac {1} { 1 + K |r|   } \cdot \frac {1} { 1 + P |\sin \omega | } \rp^{A}   
\end{align*}
for any $i, j, A \geqslant 0$.   Further, we may truncate the integral at $ \pm (K^{\vepsilon}/K, P^{\vepsilon}/P) $ at the cost of an exponentially small error $O (|z|^2 (K+P)^{-A} )$.   }

\begin{lem}\label{cor: 1}
We have
\begin{align}\label{4eq: integral H(z), 2}
	\SDH (z) =     |2  z|^2   \viint_{ \SHA  }  \cos (2 \Re (z \trh (r, \omega)) )  (\sinh^2 r  +  \sin^2 \omega  ) k  (  r) \theta (  \omega)   \nd \widehat{\mu} (r,  \omega ) . 
\end{align} 
\end{lem}



\section{Proof of Theorem \ref{thm: large sieve}} 

We start by multiplying both sides of  the Kuznetsov formula \eqref{1eq: Kuznetsov}  by ${a}_{(m)} \overline{a}_{(n)}$ and then summing over $m, n \in \CaloO$ with $  {N} < |m|, |n| \leqslant \sqrt{2}N $.  By the choice of $ h (\vkappa, p) $ as in \eqref{4eq: choice of h},  it follows that 
\begin{align}\label{5eq: first bound}
\SC (K, P, N)  + \SE (K, P, N)  	\Lt K P \big(K^2+P^2\big) \|\boldsymbol{a}_N  \|_2^2 + \sum_{c } \frac {|\varPhi (c, K, P, N) |}  {|c|^2}  , 
\end{align}
with
\begin{align*}
\varPhi (c, K, P, N) =	\mathop{\sum\sum}_{m, n}   {a}_{(m)} \overline{a}_{(n)} S (m, n;c) \SDH \bigg( \frac {2\pi \sqrt{mn}} {c} \bigg) . 
\end{align*}

Next, we apply Proposition \ref{prop: quad form} and Lemmas \ref{prop: integral repn}, \ref{cor: 1} to estimate $\varPhi (c, K, P, N) $ in various ranges of $c$.   Note that $|\theta| = |\trh (r, \omega)| = \sqrt{\sinh^2 r + \cos^2 \omega}$.

If $|c| \geqslant N^2$ then by  \eqref{4eq: integral H(z), 2} and \eqref{1eq: Kloosterman, 1} we get
\begin{align*}
	 \varPhi (c, K, P, N)  \Lt \lp \frac 1 {K^2} + \frac 1 {P^2} \rp  \frac {\tau^2 (c) N^4} {|c| }   \|\boldsymbol{a}_N \|_2^2. 
\end{align*}

If $ {N} \leqslant |c| < N^2$ then by  \eqref{4eq: integral H(z), 2} and \eqref{1eq: Kloosterman, 2} we get
\begin{align*}
	\varPhi (c, K, P, N)  \Lt \lp \frac 1 {K^2} + \frac 1 {P^2} \rp  N^2 \log (3N) \|\boldsymbol{a}_N \|_2^2. 
\end{align*}

In the case $ |c| <  {N}$, we use \eqref{1eq: Kloosterman, 3} if $|r| < 1/3$ and $| \omega| < 1$,  and \eqref{1eq: Kloosterman, 2} if otherwise. 

If $ {N} / KP \leqslant |c| <  {N}$ then by   \eqref{4eq: integral H(z), 2},  \eqref{1eq: Kloosterman, 2}, and \eqref{1eq: Kloosterman, 3}  we get
\begin{align*}
	\varPhi (c, K, P, N)  \Lt \left( \bigg(  \frac {e^{ - K^2/9 }} {K} + \frac { e^{-  P^2}} {P} \bigg) \frac {N^{4 }  \log (3N)} {\, |c|^2  } + \lp \frac 1 {K^2} + \frac 1 {P^2} \rp  \frac {N^{3+ \vepsilon}} {|c| }   \right)   \|\boldsymbol{a}_N \|_2^2. 
\end{align*}

If $|c| <  {N} / KP  $ then by   \eqref{4eq: integral H(z)}, \eqref{1eq: Kloosterman, 2}, and \eqref{1eq: Kloosterman, 3}   we get
\begin{align*}
\varPhi (c, K, P, N)  \Lt  \big(K^2+P^2\big)  \lp    \big(K   e^{ - K^2/9 } + P e^{-  P^2} \big) N^2 \log (3N)	+    |c| {N^{1 + \vepsilon}}     \rp   \|\boldsymbol{a}_N \|_2^2 . 
\end{align*}

Gathering the above estimates together, we obtain 
\begin{align*}
	\sum_{c  } \frac {|\varPhi (c, K, P, N) |}  {|c|^2} \Lt \big(K^2+P^2\big) N^{2+\vepsilon} \lp  K e^{ - K^2/9 } +   P e^{-  P^2}  + \frac 1 {KP}  \rp    \|\boldsymbol{a}_N \|_2^2 ,
\end{align*}
if $ {N} \geqslant KP$, and
\begin{align*}
	\sum_{c  } \frac {|\varPhi (c, K, P, N) |}  {|c|^2} \Lt N^{2+\vepsilon} \lp    \bigg( \frac{e^{ - K^2/9 }} {K}+ \frac {e^{-  P^2}} {P}  \bigg)  N^2 +  \lp \frac 1 {K^2} + \frac 1 {P^2} \rp  {N}   \rp   \|\boldsymbol{a}_N \|_2^2 ,
\end{align*}
if $  {N} < KP$. Hence by \eqref{5eq: first bound} 
\begin{align*}
	   (\SC + \SE) (K, P, N)    	\Lt  \big(K^2+P^2\big) \bigg( \hskip -1pt  K P  +  \frac{N^{2+\vepsilon} } {KP} \big( 1 + K^2 P e^{-K^2/9} + KP^2e^{-P^2} \big) \hskip -1pt  \bigg) \|\boldsymbol{a}_N  \hskip -1pt  \|_2^2 .
\end{align*}
Since the left-hand side is non-decreasing in $K$ and $P$, we may replace $K$ or $P$ in the right-hand side by any   $K_1>K$ or $P_1 > P$. On putting $K_1 = K + P^{\vepsilon}$ or $P_1 = P + K^{\vepsilon}$, it yields
\begin{align*}
	\SC (K, P, N)    +   \SE (K, P, N)  \Lt    \big(K^2+P^2\big) \bigg(    K P + K^{1+\vepsilon} + P^{1+\vepsilon}   +  \frac{N^{2+\vepsilon} } {KP}     \bigg)   \|\boldsymbol{a}_N  \hskip -1pt  \|_2^2 .
\end{align*}
This completes the proof of Theorem \ref{thm: large sieve}.

\section{Proof of Theorem \ref{thm: moments}}\label{sec: Proof Theorem 3}

For simplicity, let us focus on the cusp-form contribution to $ \SM_{2q} (K, P) $ ($q=1, 2$). 
Let $K, P \Gt 1$ be large. By `negligibly small' we mean $ O_A ( (K^2+P^2)^{-A} )$ for any $A > 0$. 

Let $\uppii \in \Pi_{c} $.  For $\Re (s) > 1$, we have 
\begin{align*}
	L (s, \uppii) = \sum_{ (n) \subset \CaloO}  \frac {\lambda_{\uppii} (  n)  } {|n|^{2  s} }, \qquad  {L (s,  \uppii   )^2}  = {\zeta (2s, 0) } \sum_{ (n) \subset \CaloO} \frac {\lambda_{\uppii} (n) \tau (n) } {|n|^{2 s}};
\end{align*}
the latter (note that $\zeta (s, 0)$ is the Dedekind $\zeta$ function) is due to the Hecke relation 
\begin{align*}
	\lambda_{\uppii} ( m) \lambda_{\uppii}  (n) = \sum_{ (d) | (m, n) } \lambda_{\uppii}  (m n / d^2  ). 
\end{align*} 
By \cite[Theorem 5.3]{IK}, we have the following approximate functional equations:  
\begin{align*}
	L  (  1 / 2, \uppii  )    =  ( & 1+(-1)^{p_{\uppii}}) \sum_{ (n) \subset \CaloO}  \frac {   \lambda_{\uppii}  (n )  } { |n|  }     V_1   ( |  \pi n |^2 ; \vkappa_{\uppii}, p_{\uppii} )    ,  \\
 	L  (  1 / 2,   \uppii   )^2 &   =  2 \sum_{ (n) \subset \CaloO }  \frac {   \lambda_{\uppii}  (n ) \tau (n) } { |n|  }     V_2     ( |\pi^2 n |^2 ; \vkappa_{\uppii}, p_{\uppii}   )  , 
\end{align*}
in which $(-1)^{p_{\uppii}}$ is the root number of $\uppii$. 

Now let   $\uppii \in \Pi_c (K, P)$. The $\gamma$ factor of $\uppii$ is equal to $(2\pi)^{-2s} \gamma (s, \vkappa_{\uppii}, p_{\uppii})$, with 
\begin{align*}
	\gamma (s, \vkappa, p) = \Gamma \lp s + \frac {i\vkappa +|p|} 2 \rp \Gamma  \lp s + \frac {- i\vkappa +|p|} 2 \rp. 
\end{align*} Thus  \cite[Proposition 5.4]{IK} implies that  one may effectively restrict the sums above to the ranges $ |n| \leqslant (K^2+P^2)^{q/2+\vepsilon} $ ($q=1, 2$) at the cost of a negligibly small error. 

A minor difficulty arises for the application of Theorem \ref{thm: large sieve} as $V_q ( |\pi^q n |^2 ; \vkappa_{\uppii}, p_{\uppii} )$ depends on $(\vkappa_{\uppii}, p_{\uppii})$. To address this issue, we invoke the following expressions due to Blomer \cite[Lemma 1]{Blomer}:
\begin{align*}
	V_1  (y; \vkappa,  p)  =   \frac 1 {2   \pi i }   \int_{ \vepsilon - i U}^{\vepsilon + i U}  & G  (v, \vkappa,  p)    y^{ - v}   \frac {\nd v} {v} +  {O_{\vepsilon }} \bigg( \frac { (K^2+P^2)^{  \vepsilon} } {y^{ \vepsilon} \exp ({U^2 / 2}) } \bigg),  \\
	V_2  (y; \vkappa,  p)  = \frac 1 {2   \pi i }   \int_{ \vepsilon - i U}^{\vepsilon + i U}     \zeta (1+&2v, 0) G  (v, \vkappa,  p)^2     y^{ - v}   \frac {\nd v} {v} +  {O_{\vepsilon}} \bigg( \frac { (K^2+P^2)^{  \vepsilon} } {y^{ \vepsilon} \exp ({U^2  }) } \bigg) ,
\end{align*}  
where 
\begin{align*}
	G (v, \vkappa,  p) = 
	\frac {\gamma (1/2+ v, \vkappa, p)} {\gamma (1/2, \vkappa, p)} \exp({v^2})   . 
\end{align*}  
 The error terms above are again negligibly small if we choose $ U = \log ( K^2+P^2)$. By the Stirling formula, for any $v$ on the integral contour and $ (\vkappa, p) $ in the rectangle $  \BSA (K, P) $ 
 \begin{align*}
 	G  (v, \vkappa,  p) = O_{\vepsilon} \big( (K^2+P^2)^{\vepsilon}\big). 
 \end{align*} 
 
 By a dyadic partition of summation, Cauchy--Schwarz, and an application of Theorem \ref{thm: large sieve} inside the integral, one may easily derive the bounds
 \begin{align*}
 &	\sum_{\uppii \shskip \in \Pi_{c} (K, P) }  \frac { L (1/2,  \uppii)^2  } { L (1, \mathrm{Sym}^2 \uppii) }   \Lt   K P \big(K^2+P^2\big)^{1+\vepsilon}, \\ & \quad \sum_{\uppii \shskip \in \Pi_{c} (K, P) }  \frac { L (1/2,  \uppii)^4  } { L (1, \mathrm{Sym}^2 \uppii) }   \Lt \frac {(K^2+P^2)^{3+\vepsilon}} {KP}. 
 \end{align*}
For the second inequality, since the left-hand side is non-decreasing in $K$ and $P$, we may replace $K$ and $P$ in the right-hand side by $K+P$ so that it is improved into
\begin{align*}
	\sum_{\uppii \shskip \in \Pi_{c} (K, P) }  \frac { L (1/2,  \uppii)^4  } { L (1, \mathrm{Sym}^2 \uppii) }   \Lt   {\big(K^2+P^2 \big)^{2+\vepsilon}}  . 
\end{align*}

\delete{Similarly,
\begin{align*} 
 &	\zeta   (s+i \vkappa,  p )   \zeta   (s-i \vkappa,  p ) = \sum_{ (n) }  \frac {\tau_{i \vkappa,  p} (  n)  } {|n|^{2  s} } , \\
&	\zeta   (s+i \vkappa,  p)^2   \zeta    (s-i\vkappa,  p )^2 = \zeta (2s,0) \sum_{ (n)}  \frac {\tau_{i \vkappa,  p} (  n)  \tau (n)} {|n|^{2  s} } .
\end{align*}
}


\def\cprime{$'$}

\end{document}